\documentclass{article}

\usepackage[cp1251]{inputenc}
\usepackage[russian]{babel}
\usepackage{amsmath}
\usepackage{amsfonts}
\usepackage{latexsym}
\usepackage{graphicx}
\usepackage{amsthm}

\def\Char{\operatorname{char}}

\def\A{\operatorname{A}}
\def\D{\operatorname{D}}
\def\E{\operatorname{E}}
\def\GL{\operatorname{GL}}

\def\a{\alpha}
\def\b{\beta}
\def\g{\gamma}
\def\d{\delta}

\def\l{\lambda}

\def\Z{{\Bbb Z}}

\newtheorem{lemma}{Лемма}
\newtheorem{theorem}{Теорема}
\newtheorem*{theorem*}{Теорема}

\begin{document}

\begin{center}
{\bf  Орбиты векторов некоторых представлений. III}

\vspace{2mm}

И.М. Певзнер \footnote[1]{Настоящая работа выполнена при
содействии проекта РФФИ 19-01-00297}

\end{center}

\vspace{1cm}

Пусть $\Phi$~--- система корней одинаковой длины и $K$~---
произвольное поле. Далее, обозначим через $\d$  максимальный
корень $\Phi$ и положим $\Phi_0 = \{\a\in\Phi; \d\perp\a\}$, $G_0
= G_{\text{\rm sc}}(\Phi_0, K)$ и $V_1 = \langle e_\a; \angle(\a,
\d) = \pi/3\rangle$, где $e_\a$~--- элементарные корневые
элементы. В настоящей серии статей рассмотрены орбиты
действия $G_0$ на $V_1$.

Такое действие изучалось во множестве работ. Прежде всего это,
разумеется, случай $\Phi=\E_8$~--- тогда получается $56$-мерное
минимальное микровесовое представление группы типа $\E_7$.
Остальные случаи исследуются меньше, однако тоже встречаются
достаточно часто. Наряду с изучением собственно представления
$G_0$ в $V_1$, это может помочь и при исследовании представления
всей группы $G_{\text{\rm sc}}(\Phi, K)$ и соответствующей алгебры
Ли.

Настоящая статья является продолжением работ [18, 19]. В них были доказаны некоторые общие результаты и
разобраны случаи $\Phi=\E_l, \A_l$ и $\D_l$ при $\Char K\neq 2$. В настоящей работе будут
рассмотрены случаи $\Phi=\E_l$ при $\Char K = 2$.

\begin{center}
{\bf \S1. Основные обозначения}
\end{center}

Пусть $K$~--- произвольное поле с характеристикой, равной $2$, $\Phi$~--- система корней одной длины, $V = V(\Phi)$~--- соответствующая алгебра Ли,  а $G=G_{\text{\rm sc}}(\Phi, K)$~--- соответствующая односвязная группа.

Как известно, в $V$ существует базис Шевалле $\{e_\a,\a\in\Phi;
h_\a,\a\in\Pi\}$, где $\Pi$~-- фундаментальная система корней. При
этом все $h_\a$ из подалгебры Картана; $[e_\a,e_{-\a}] = h_\a$;
$[h_\a,e_\b] = A_{\a\b}e_\b$, где $A_{\a\b} =
2(\a,\b)/(\a,\a)\in\Z$~--- числа Картана; $[e_\a,e_\b] =
N_{\a\b}e_{\a+\b}$ при $\a+\b\in\Phi$ и $[e_\a,e_\b]=0$ при
$\a+\b\notin\Phi$ и $\b\neq -\a$, где $N_{\a\b} = \pm1$~---
структурные константы. В рассматриваемом случае $\Char K=2$ они все равны $1$, но в этом параграфе мы их все же оставим для полноты картины. Коэффициент в разложении вектора $x\in V$ по этому базису при $e_\a$ обозначим $x^\a$, а соответствующий
элемент из подалгебры Картана обозначим $x^h$; тогда $x =
\sum_{\a\in\Phi} x^\a e_\a + x^h$. 

Далее, в группе $G$ выделяются элементарные корневые элементы
$x_\a(a)$, $\a\in\Phi, a\in K$ и $X_\a = \langle x_\a(a); a\in
K\rangle$~--- элементарные корневые подгруппы. В работе будут
использоваться формулы для действия $x_\a(a)$ на базисе Шевалле.
Они перечислены, например, в [35]. Нам понадобятся следующие
равенства: $x_\a (a) e_\b = e_\b$ при $\angle(\a,\b)<2\pi/3$,
$x_\a (a) e_\b = e_\b + N_{\a\b}  a e_{\a+\b}$ при
$\angle(\a,\b)=2\pi/3$, $x_\a (a) e_{-\a} = e_{-\a} + a h_\a - a^2
e_\a$ и $x_\a (a) h_\b = h_\b - A_{\b\a} a e_\a$; мы будем ими
пользоваться без дополнительных ссылок.

Обозначим через $\d$ максимальный корень системы $\Phi$.
Экстраспециальным унипотентным радикалом называется подгруппа
$U_\d = \langle X_\a; \angle(\a,\d)<\pi/2\rangle$. Подробнее об
этом радикале говорится, например, в [7, 16].

Пусть $\a\in\Phi$~--- некоторый корень. Разобьем все корни из
$\Phi$ на пять классов в зависимости от их расположения
относительно корня $\a$: $\Phi_2(\a) = \{\a\}$, $\Phi_1(\a) =
\{\b; \angle(\b,\a) = \pi/3\}$, $\Phi_0(\a) = \{\b; \angle(\b,\a)
= \pi/2\}$, $\Phi_{-1}(\a) = \{\b; \angle(\b,\a) = 2\pi/3\}$,
$\Phi_{-2}(\a) = \{-\a\}$. Другими словами, $\b$ принадлежит
$\Phi_i(\a)$ тогда и только тогда, когда скалярное произведение
$\b$ и $\a$ равно $i/2$. Простейшие свойства этого разбиения
приведены в [16]. Часто нас будет интересовать случай $\a=\d$; для
краткости, аргумент у $\Phi_i(\d)$ будем опускать. В этих
обозначениях экстраспециальный радикал $U_{\d} = \langle X_{\a};
\a\in\Phi_2\cup\Phi_1\rangle$.

\begin{center}
{\bf \S2. Необходимые результаты из работы [18]}
\end{center}

1. Так как настоящая работа является переносом результатов статьи
[18] со случая $\Char K\neq 2$ на случай $\Char K=2$,
стоит напомнить основные моменты той статьи. Мы рассматриваем
орбиты действия группы $G_0 = G_{\text{\rm sc}}(\Phi_0, K)$ на
пространстве $V_1 = \langle e_\a; \a\in\Phi_1\rangle$. Пусть $x\in
V_1$. В начале мы доказали, что можно считать, что $x=x^\l e_\l +
x^{\d-\l} e_{\d-\l} + x^\mu e_\mu + x^\nu e_\nu + x^\xi e_\xi$;
здесь $\l, \mu, \nu, \xi \in \Phi_1$~--- некоторые попарно
ортогональные корни и $\l+\mu+\nu+\xi = 2\d$. Далее мы доказывали,
что для систем $\Phi=\E_6$, $\E_7$ и $\E_8$ корни $\l$, $\mu$, $\nu$
и $\xi$ могут быть выбраны произвольно. После
этого мы заметили, что любой $x\in V_1$ приводится к одному из
следующих типов.

\begin{enumerate}
\item[I.] $x=0$.

\item[II.] $x=x^\l e_\l$; $x^\l\neq 0$. Есть для всех систем
$\Phi$, кроме $\A_1$.

\item[III.] $x=x^\l e_\l + x^\mu e_\mu$; $x^\l, x^\mu\neq 0$. Есть
для всех систем $\Phi$, кроме $\A_1$ и $\A_2$.

\item[IV.] $x=x^\l e_\l + x^\mu e_\mu + x^\nu e_\nu$; $x^\l,
x^\mu, x^\nu\neq 0$. Есть для всех систем $\Phi$, кроме $\A_l$.

\item[V.] $x=x^\l e_\l + x^\mu e_\mu + x^\nu e_\nu +x^\xi e_\xi$;
$x^\l, x^\mu, x^\nu, x^\xi\neq 0$. Есть для всех систем $\Phi$,
кроме $\A_l$.

\item[VI.] $x=x^\l e_\l + x^{\d-\l} e_{\d-\l}$; $x^\l,
x^{\d-\l}\neq 0$. Есть для всех систем $\Phi$, кроме $\A_1$.

\item[VII.] $x=x^\l e_\l + x^{\d-\l} e_{\d-\l} + x^\mu e_\mu$;
$x^\l, x^{\d-\l}, x^\mu\neq 0$. Есть для всех систем $\Phi$, кроме
$\A_1$ и $\A_2$.

\item[VIII.] $x=x^\l e_\l + x^{\d-\l} e_{\d-\l} + x^\mu e_\mu +
x^\nu e_\nu$; $x^\l, x^{\d-\l}, x^\mu, x^\nu\neq 0$. Есть для всех
систем $\Phi$, кроме $\A_l$.

\item[IX.] $x=x^\l e_\l + x^{\d-\l} e_{\d-\l} + x^\mu e_\mu +
x^\nu e_\nu +x^\xi e_\xi$; $x^\l, x^{\d-\l}, x^\mu,
x^\nu,x^\xi\neq 0$. Есть для всех систем $\Phi$, кроме $\A_l$.
\end{enumerate}

От VII и VIII случаев мы избавились сразу же. Далее, от VI случая
мы избавлялись при $\Phi\neq \A_l$ и $|K|>2$, а от IX~--- при
$\Char K\neq 2$. Кроме того, мы доказали следующую несложную
лемму:

\begin{lemma}[Лемма 1 из {[18]}]
Пусть
$w_{\a}(a)=x_{-\a}(-a^2+a)x_{\a}(-\frac1a)x_{-\a}(a-1)x_{\a}(1)$
при $\a\in\Phi_0$ и $a\in K^*$. Тогда $w_\a(a)e_\b = \frac1a e_\b$
при $\angle(\a,\b) = \pi/3$; $w_\a(a)e_\b = e_\b$ при
$\angle(\a,\b) = \pi/2$ и $w_\a(a)e_\b = a e_\b$ при
$\angle(\a,\b) = 2\pi/3$.
\end{lemma}

С помощью этой леммы мы избавлялись от большинства коэффициентов в
оставшихся случаях и приходили к искомому списку орбит.

2. Далее мы доказывали, что все описанные случаи дают разные
орбиты. Во-первых, была доказана следующая лемма:

\begin{lemma}[Лемма 3 из {[18]}]
Пусть $\Phi\neq \A_l$, $\Char K\neq 2$ и $x = \sum_{\a\in\Phi_1}
x^\a e_\a$. Тогда существует и единственен корневой элемент $y =
\sum_{\a\in\Phi} y^\a e_\a + y^h$, в котором $y^\a=x^\a$ при
$\a\in\Phi_1$, $y^\d = 1$ и $y^h\in\langle h_\a;
\a\in\Phi_0\cap\Pi\rangle$.
\end{lemma}

Во-вторых, мы вводили следующие определения: элемент $x$
называется темным, если угол между соответствующим ему $y$ и
$e_\d$ равен $\pi$ (иначе говоря, $y^{-\d}\neq 0$). Далее, $x$
называется светящимся, если этот угол равен $2\pi/3$; блестящим,
если угол равен $\pi/2$ и сингулярным, если он равен $\pi/3$. Как
несложно видеть, последнее условие равносильно тому, что $x$
корневой элемент. Названия взяты из работы [30] для случая $\Phi =
\E_8$. Их определения в [30] другие, но, на самом деле,
равносильные.

Затем мы заметили, что при умножении на $g\in G_0$ элементу $gx$
соответствует корневой элемент $gy$. При этом угол между $y$ и
$e_\d$ при таком умножении также не меняется, значит определения
темного, светящегося, блестящего и сингулярного элементов можно
расширить на орбиты. Более того, если $x$ темный, то коэффициент
$y^{-\d}$ также не меняется при умножении на $g\in G_0$, то есть
тоже является инвариантом. Наконец, мы проверяли, что в случае II
вектор $x$ получается сингулярным, в случае III~--- блестящим, в
случае IV~--- светящимся, а в случае V~--- темным. При этом в
случае V инвариант $y^{-\d}$ равнялся произведению нескольких
структурных констант и $x^\l x^\mu x^\nu x^\xi$, поэтому
произведение всех коэффициентов $x^\l x^\mu x^\nu x^\xi$ также
постоянно. Для случая $\Char K\neq 2$, рассматриваемого в статье [18], этого
хватало для классификации. К сожалению, не все вышеописанное рассуждение проходит для рассматриваемого в настоящей работе случая $\Char K=2$, хотя оно по-прежнему остается очень полезным.

\begin{center}
{\bf \S3. Конструкция орбит при $\Char K=2$ и $|K|>2$}
\end{center}	

3. Согласно 1 пункту, в случае $\Phi=\E_l$, $\Char K=2$ и $|K|>2$ любой $x\in V_1$ приводится к одному из
следующих типов.

\begin{enumerate}
\item[I.] $x=0$.

\item[II.] $x=x^\l e_\l$; $x^\l\neq 0$.

\item[III.] $x=x^\l e_\l + x^\mu e_\mu$; $x^\l, x^\mu\neq 0$.

\item[IV.] $x=x^\l e_\l + x^\mu e_\mu + x^\nu e_\nu$; $x^\l,
x^\mu, x^\nu\neq 0$.

\item[V.] $x=x^\l e_\l + x^\mu e_\mu + x^\nu e_\nu +x^\xi e_\xi$;
$x^\l, x^\mu, x^\nu, x^\xi\neq 0$.

\item[IX.] $x=x^\l e_\l + x^{\d-\l} e_{\d-\l} + x^\mu e_\mu +
x^\nu e_\nu +x^\xi e_\xi$; $x^\l, x^{\d-\l}, x^\mu,
x^\nu,x^\xi\neq 0$.
\end{enumerate}

Теперь посмотрим, как можно упростить эти типы с помощью леммы 1. Напомним также, что корни $\l, \mu, \nu, \xi \in \Phi_1$ мы можем выбирать сами (с учетом попарной ортогональности), что упрощает подбор корня $\a$ из этой леммы. В типе I, разумеется, остается $0$. В типе II, выбирая корень $\a$ из леммы так, чтобы $\angle(\l,\a)=\pi/3$, можно сделать $x^\l$ равным $1$. В типе III ситуация схожая~--- выбирая $\a$ так, чтобы $\angle(\l,\a)=\pi/3$ и $\angle(\mu,\a)=\pi/2$, можно сделать $x^\l$ равным $1$, а затем, выбирая $\a$ так, чтобы $\angle(\l,\a)=\pi/2$ и $\angle(\mu,\a)=\pi/3$, можно сделать и $x^\mu$ равным $1$. Действуя аналогично, в типе IV можно сделать $x^\l = x^\mu = x^\nu = 1$. В типе V ситуация немного сложнее, так как не существует корня $\a\in\Phi_0$, ортогонального к $\mu, \nu, \xi$ и неортогонального к $\l$. Выбирая корень $\a$ так, чтобы $\angle(\l,\a)=\pi/3$, $\angle(\mu,\a)=\angle(\nu,\a)=\pi/2$ и $\angle(\xi,\a)=2\pi/3$, можно сделать $x^\l$ равным $1$; при этом $x^\mu$ и $x^\nu$ не меняются. Аналогичными действиями можно сделать $x^\mu$ и $x^\nu$ равными $1$, при этом $x^\xi$ станет равным произведению всех четырех изначальных коэффициентов. Повторяя рассуждения для типа IX, можно и в этом случае сделать коэффициенты $x^\l$, $x^\mu$ и $x^\nu$ равными $1$; при этом коэффициент $x^{\d-\l}$ станет равен произведению изначальных коэффициентов при корнях $\l$ и $\d-\l$, а коэффициент $x^\xi$~--- произведению изначальных коэффициентов при корнях $\l$, $\mu$, $\nu$ и $\xi$.

4. Типы V и IX, оказывается, можно и еще упростить. Для этого стоит вспомнить, как в статье [18] для случая $\Char K\neq 2$ избавлялись от типа IX. Повторим это рассуждение, слегка поменяв обозначения под наши цели. В той работе было показано, что можно считать все используемые в рассуждении структурные константы равными $1$. В рассматриваемом сейчас случае все еще проще, при $\Char K = 2$ структурные константы не важны, поэтому мы их не записываем.

Пусть $x_1=x=x^\l e_\l + x^{\d-\l} e_{\d-\l} + x^\mu e_\mu + x^\nu e_\nu + x^\xi e_\xi$; $x^\l, x^\mu, x^\nu,x^\xi\neq 0$. ``Испортим''\ этот элемент, а затем снова приведем его к такому же виду. Пусть $k\in K$ и $x_2 = x_{\d-\mu-\nu}(k) x_1 = x^\l e_\l + x^{\d-\l} e_{\d-\l} + x^\mu e_\mu + x^\nu e_\nu + (x^\xi + k x^{\d-\l}) e_\xi + k x^\mu e_{\d-\nu} + k x^\nu e_{\d-\mu}$. Теперь по очереди избавляемся в этом выражении от $e_{\d-\nu}$ и $e_{\d-\mu}$. Пусть $x_3 = x_{\d-\l-\nu}(k \frac{x^\mu}{x^\l}) x_2 = x^\l e_\l + (x^{\d-\l} + k \frac{x^\mu x^\nu}{x^\l}) e_{\d-\l} + x^\mu e_\mu + x^\nu e_\nu + (x^\xi + k x^{\d-\l} + k^2 \frac{x^\mu x^\nu}{x^\l}) e_\xi + k x^\nu e_{\d-\mu}$. Наконец, пусть $x_4 = x_{\d-\l-\mu}(k \frac{x^\nu}{x^\l}) x_3 = x^\l e_\l + x^{\d-\l} e_{\d-\l} + x^\mu e_\mu + x^\nu e_\nu + (x^\xi + k x^{\d-\l} + k^2 \frac{x^\mu x^\nu}{x^\l}) e_\xi$.

В рассуждениях из пункта 3  у нас получалось в этих типах два ``инварианта''~--- произведение изначальных коэффициентов при корнях $\l$ и $\d-\l$ и произведение изначальных коэффициентов при корнях $\l$, $\mu$, $\nu$ и $\xi$. При преобразованиях из предыдущего абзаца первое произведение остается прежним, а второе меняется. Запишем это изменение так, чтобы оно не зависело от коэффициентов разложения нашего элемента, а только от первого ``инварианта''\ (как мы увидим позднее, это действительно инвариант). Произведение $x^\l x^\mu x^\nu x^\xi$ перешло в $x^\l x^\mu x^\nu (x^\xi + k x^{\d-\l} + k^2 \frac{x^\mu x^\nu}{x^\l}) = x^\l x^\mu x^\nu x^\xi + x^\l x^\mu x^\nu k x^{\d-\l} + x^\mu x^\nu k^2 x^\mu x^\nu$. Полагая $l = x^\mu x^\nu k$, получаем $x^\l x^\mu x^\nu x^\xi + x^\l x^{\d-\l} l + l^2$.

5. Для фиксированного $s\in K$ рассмотрим отношение $a\sim_s b\Leftrightarrow \exists l\in K : a = b+l s+l^2$. Несложно видеть, что для рассматриваемого случая $\Char K = 2$ это будет отношение эквивалентности. Обозначим через $K_s$ множество представителей классов эквивалентности по этому отношению.

Тогда при $\Phi=\E_6, \E_7$ или $\E_8$, $\Char K=2$ и $|K|>2$ любой $x\in V_1$ приводится к одному из следующих типов:

\begin{enumerate}
\item[I.] $x=0$.

\item[II.] $x=e_\l$.

\item[III.] $x=e_\l + e_\mu$.

\item[IV.] $x=e_\l + e_\mu + e_\nu$.

\item[V.] $x=e_\l + e_\mu + e_\nu +x^\xi e_\xi$;
$x^\xi\in K_0$, $x^\xi\notin \bar 0$.

\item[IX.] $x=e_\l + s e_{\d-\l} + e_\mu +
e_\nu + x^\xi e_\xi$; $s\in K^*, x^\xi\in K_s$.
\end{enumerate}

\begin{center}
{\bf \S4. Доказательство различности орбит при $\Char K=2$ и $|K|>2$}
\end{center}	

6. Осталось доказать, что все полученные типы лежат в разных орбитах. Для этого повторим, с соответствующими изменениями, доказательство леммы~3 из [18]. 

Как говорилось в утверждении 2 [17], любой корневой элемент $y = \sum_{\a\in\Phi} y^\a e_\a + y^h$ с $y^\d = 1$ представляется в виде $y = ue_\d$, где $u\in U_{-\d}$~--- унипотентный элемент, равный $$u = x_{-\d}(a) \cdot\prod\limits_{\g\in\Phi_{-1}} x_\g\Big(N_{\g\d} y^{\d+\g}\Big)$$ для некоторого $a\in K$. Таким образом, для любого $x\in V_1$ существует целая серия корневых элементов $y$, таких что $y^\d = 1$ и $y^\a=x^\a$ при $\a\in\Phi_1$, элементы которой отличаются друг от друга выбором коэффициента $a$.  
Есть ровно один простой корень $\a_k$, не ортогональный к $\d$ (для $\Phi=\E_6$ это $\a_2$, для $\Phi=\E_7$~--- $\a_1$, а для $\Phi=\E_8$~--- $\a_8$), причем в разложении $\d$ этот корень входит с коэффициентом $2$. Соответственно, разложим $y^h$ по базису $h_i$ и посмотрим на коэффициент при $h_k$. В лемме~3 [18] утверждалось, что при $\Char K\neq 2$ можно этот коэффициент сделать нулевым; в нашем случае ситуация немного другая.

7. Пусть $y_0$~--- это один из корневых элементов, таких что $y_0^\d = 1$ и $y_0^\a=x^\a$ при $\a\in\Phi_1$, а $y_1 = x_{-\d}(a)y_0$. Тогда, по описанным в \S1 формулам, получаем, что $y_1^h = y_0^h + a h_{-\d}$. Так как в разложении $\d$ корень $\a_k$ входит с коэффициентом $2$, то коэффициент при $h_k$ у элемента $y_1^h$ отличается от соответствующего коэффициента $y_2^h$ на $-2a$. Поэтому при $\Char K=2$ коэффициент элемента $y$ при $h_k$ не зависит от выбора $a$, то есть определяется элементом $x\in V_1$ однозначно. Обозначим этот коэффициент буквой $t$.

Далее, посмотрим, как меняется коэффициент $y^{-\d}$. Пусть снова $y_0$~--- это один из корневых элементов, таких что $y_0^\d = 1$ и $y_0^\a=x^\a$ при $\a\in\Phi_1$, а $y_1 = x_{-\d}(a)y_0$. Тогда $y_1^{-\d} = y_0^{-\d} - a^2 + a t$. Таким образом, элементу $x$ можно сопоставить число $t\in K$, являющееся коэффициентом соответствующего элемента $y$ при $h_k$, и класс эквивалентности по отношению  $\sim_t$ в обозначениях пункта 5. Несложно видеть также, что оба параметра не меняются при действии группой $G_0$, как и в [18]. Поэтому элементы $x$ из одной орбиты имеют одинаковые значения параметров.

8. Свяжем эти параметры с коэффициентами из 5 пункта. Пусть $x=x^\l e_\l + x^{\d-\l} e_{\d-\l} + x^\mu e_\mu + x^\nu e_\nu +x^\xi e_\xi$. Так как в рассматриваемом случае $\Char K=2$, то структурные константы можно сразу опустить; для краткости также уберем и $x_{-\d}(a)$, то есть положим $a$ равным $0$. Тогда 

\begin{multline}
y_0 = u e_\d = x_{-\l}(x^{\d-\l}) x_{\xi-\d}(x^\xi) x_{\nu-\d}(x^\nu) x_{\mu-\d}(x^\mu) x_{\l-\d}(x^\l) e_\d = \\
= x_{-\l}(x^{\d-\l}) x_{\xi-\d}(x^\xi) x_{\nu-\d}(x^\nu) x_{\mu-\d}(x^\mu) (e_\d + x^\l e_\l) = \\
= x_{-\l}(x^{\d-\l}) x_{\xi-\d}(x^\xi) x_{\nu-\d}(x^\nu) (e_\d + x^\l e_\l + x^\mu e_\mu + x^\mu x^\l e_{\mu+\l-\d}) = \\
= x_{-\l}(x^{\d-\l}) x_{\xi-\d}(x^\xi) (e_\d + x^\l e_\l + x^\mu e_\mu + x^\mu x^\l e_{\mu+\l-\d} + x^\nu e_\nu + x^\nu x^\l e_{\nu+\l-\d}+ \\
+ x^\nu x^\mu e_{\nu+\mu-\d} + x^\nu x^\mu x^\l e_{-\xi}) = x_{-\l}(x^{\d-\l}) (e_\d + x^\l e_\l + x^\mu e_\mu + x^\mu x^\l e_{\mu+\l-\d} + \\
+ x^\nu e_\nu + x^\nu x^\l e_{\nu+\l-\d}+ x^\nu x^\mu e_{\nu+\mu-\d} + x^\nu x^\mu x^\l e_{-\xi} + x^\xi e_\xi + x^\xi x^\l e_{\xi+\l-\d} + \\
+ x^\xi x^\mu e_{\xi+\mu-\d} + x^\xi x^\mu x^\l e_{-\nu} + x^\xi x^\nu e_{\xi+\nu-\d} + x^\xi x^\nu x^\l e_{-\mu} + x^\nu x^\nu x^\mu e_{-\l} + x^\xi x^\nu x^\mu x^\l e_{-\d}) = \\
= e_\d + x^\l e_\l + x^\mu e_\mu + x^\mu x^\l e_{\mu+\l-\d} 
+ x^\nu e_\nu + x^\nu x^\l e_{\nu+\l-\d}+ x^\nu x^\mu e_{\nu+\mu-\d} + \\
+ x^\nu x^\mu x^\l e_{-\xi} + x^\xi e_\xi + x^\xi x^\l e_{\xi+\l-\d} + 
x^\xi x^\mu e_{\xi+\mu-\d} + x^\xi x^\mu x^\l e_{-\nu} + x^\xi x^\nu e_{\xi+\nu-\d} + \\
+ x^\xi x^\nu x^\l e_{-\mu} + (x^\nu x^\nu x^\mu + (x^{\d-\l})^2 x^\l) e_{-\l} + x^\xi x^\nu x^\mu x^\l e_{-\d} + x^{\d-\l} e_{\d-\l} + x^{\d-\l} x^\l h_{-\l} + \\
+ x^{\d-\l} x^\mu x^\l e_{\mu-\d} + x^{\d-\l} x^\nu x^\l e_{\nu-\d} + x^{\d-\l} x^\xi x^\l e_{\xi-\d} 
\end{multline}

Таким образом, $y_0^{-\d}=x^\xi x^\nu x^\mu x^\l$, а $y_0^h = x^{\d-\l}x^\l h_\l$. В разложении корня $\l$ на простые корни $\a_k$ входит с коэффициентом $1$, поэтому в разложении $y_0^h$ элемент $h_k$ будет входить с коэффициентом $x^{\d-\l}x^\l$. Как мы уже говорили, этот коэффициент не меняется при умножении $y_0$ на $x_{-\d}(a)$, а $y_0^{-\d}$  остается в том же классе эквивалентности. Для интересующих нас типов IV, V и IX $x=e_\l + s e_{\d-\l} + e_\mu + e_\nu + x^\xi e_\xi$; $s\in K, x^\xi\in K_s$, значит в разложении $y^h$ элемент $h_k$ будет входить с коэффициентом $s$, а $y^{-\d} = x^\xi$. Поэтому для различных $s\in K$ и $x^\xi\in K_s$ действительно будут получаться разные орбиты.

9. Осталось научиться ``отделять''\ II, III и IV типы (как всегда, I тип выделяется сразу), так как всем этим типам соответствуют $s=0$ и $x^\xi = 0$. Заметим, что в этом случае существует ровно один $y$, такой что $y^{-\d}=0$: в самом деле, его существование следует из подсчетов из 8 пункта, а единственность~--- из формулы $y_1^{-\d} = y_0^{-\d} - a^2 + a s$ из 7 пункта. Тогда, как и в [18], можно посмотреть на угол между этим $y$ и $e_\d$, причем при умножении на $g\in G_0$ элементу $gx$ соответствует $gy$ и угол между $y$ и $e_\d$ не меняется. Таким образом, угол является инвариантом. Посмотрим, какие углы получаются в типах II, III и IV. Вместо подсчетов, проводившихся в [18], можно сослаться на формулу (1), подставляя в нее соответствующие коэффициенты. В типе II, подставляя $x^\l = 1, x^\mu = x^\nu = x^\xi = x^{\d-\l} = 0$, получаем угол, равный $\pi/3$; такие $x$ в [18] назывались сингулярными. В типе III, подставляя $x^\l = x^\mu = 1, x^\nu = x^\xi = x^{\d-\l} = 0$, получаем угол, равный $\pi/2$; такие $x$ назывались блестящими. Наконец, в типе IV, подставляя $x^\l = x^\mu = x^\nu = 1, x^\xi = x^{\d-\l} = 0$, получаем угол, равный $2\pi/3$; такие $x$ назывались светящимися.

10. Объединяя результаты 5 и 9 пунктов, получаем теорему.

\begin{theorem} Пусть $\Phi=\E_6$, $\E_7$ или $\E_8$, $\Char K = 2$ и $|K|>2$. Для фиксированного $s\in K$ введем отношение эквивалентности $a\sim_s b\Leftrightarrow \exists l\in K : a = b+l s+l^2$ и обозначим через $K_s$ множество представителей классов эквивалентности по этому отношению. Тогда векторы из $V_1$ под действием $G_0$ образуют следующие орбиты.
\begin{enumerate}
\item Нулевая орбита, $x=0$.

\item Одна орбита из сингулярных векторов. Ее элементы можно привести к виду $x=e_\l$.

\item Одна орбита из блестящих векторов. Ее элементы можно привести к виду $x=e_\l + e_\mu$.

\item Одна орбита из светящихся векторов. Ее элементы можно привести к виду $x=e_\l + e_\mu + e_\nu$.

\item Двухпараметрическое множество орбит, элементы которых можно привести к виду $x=e_\l + s e_{\d-\l} + e_\mu + e_\nu + x^\xi e_\xi$; где $s\in K, x^\xi\in K_s$, за исключением случая $s=0, x^\xi\in\bar 0$.
\end{enumerate}
\end{theorem}

\begin{center}
{\bf \S5. Случай $|K|=2$}
\end{center}

11. Согласно 1 пункту, в случае $\Phi=\E_l$ и $|K|=2$ любой $x\in V_1$ приводится к одному из следующих типов.

\begin{enumerate}
\item[I.] $x = 0$.

\item[II.] $x = e_\l$.

\item[III.] $x = e_\l + e_\mu$.

\item[IV.] $x = e_\l + e_\mu + e_\nu$.

\item[V.] $x = e_\l + e_\mu + e_\nu + e_\xi$.

\item[VI.] $x = e_\l + e_{\d-\l}$.

\item[IX.] $x = e_\l + e_{\d-\l} + e_\mu + e_\nu + e_\xi$.
\end{enumerate}

С помощью процедуры из пункта 4 можно от типа V перейти к типу IV, поэтому тип V можно убрать. 

12. Как и в предыдущем параграфе, посмотрим на соответствующие $y$. В данном случае для каждого $x$ есть два $y$, таких что $y^\d = 1$ и $y^\a=x^\a$ при $\a\in\Phi_1$: это $y_0$, посчитанный в (1), и $y_1=x_{-\d}(1) y_0$. Как отмечалось в 7 пункте, коэффициент $t$ при $h_k$ одинаков у $y_0$ и $y_1$, а $y_1^{-\d} = y_0^{-\d}+1+t$. Поэтому если $t=1$, то число $y_1^{-\d} = y_0^{-\d}$ оказывается инвариантом, а если $t=0$, то $y_1^{-\d} \neq y_0^{-\d}$ и есть ровно один $y$ с условием $y^{-\d}=0$. Как и прежде, для нахождения коэффициентов достаточно в формулу (1) подставить нужные числа. 

Тип I, как всегда, отбрасывается сразу. В типе VI получаем $y^h = h_{-\l}$, то есть $t=1$, а $y^{-\d} = 0$, а в типе IX~--- тоже $y^h = h_{-\l}$, то есть $t=1$, но $y^{-\d} = 1$. Во всех остальных случаях $y^h=0$ и, соответственно, $t=0$. Таким образом, нам осталось научиться различать между собой II, III и IV типы. В них во всех $t=0$, поэтому есть ровно один $y$ с условием $y^{-\d}=0$. Как и в 9 пункте, можно посмотреть на углы между этим $y$ и $e_\d$: в типе II угол равен $\pi/3$, в типе III~$\pi/2$, а в типе IV~--- $2\pi/3$. Соответственно, элементы снова называются сингулярными, блестящими и светящимися. В типах VI и IX буквально тем же путем действовать не получается, так как один $y$, инваариантный под действием $G_0$, выделить не получается, однако в типе VI оба $y$ образуют угол $2\pi/3$ с $e_\d$, а в типе IX~--- оба образуют угол $\pi$. Поэтому в типе VI естественно говорить про светящиеся векторы, а в типе IX~--- про темные векторы.

13. Объединяя результаты 11 и 12 пунктов, получаем теорему.

\begin{theorem} Пусть $\Phi=\E_6$, $\E_7$ или $\E_8$ и $|K|=2$. Тогда векторы из $V_1$ под действием $G_0$ образуют следующие орбиты.

\begin{enumerate}

\item Нулевая орбита, $x=0$.

\item Одна орбита из сингулярных векторов. Ее элементы можно привести к виду $x = e_\l$.

\item Одна орбита из блестящих векторов. Ее элементы можно привести к виду $x = e_\l + e_\mu$.

\item Две орбиты из светящихся векторов. Их элементы можно привести к виду $x = e_\l + e_\mu + e_\nu$ и $x = e_\l + e_{\d-\l}$.

\item Одна орбита из темных векторов. Ее элементы можно привести к виду $x = e_\l + e_{\d-\l} + e_\mu + e_\nu + e_\xi$. 
\end{enumerate}
\end{theorem}

14. Отметим напоследок, что использованные в этой работе, да и предыдущих статьях серии (хоть и в более ``закамуфлированном''\ виде), параметры~--- коэффициент соответствующего элемента $y$ при $h_k$ и $y^{-\d}$~--- не являются чем-то принципиально новым. Первый из них является квадратной формой, а второй схож с формой четвертой степени, хоть и не является ею; при этом оба параметра инвариантны под действием $G_0$. В случае $\Phi=\E_8$ действие $G_0$ на $V_1$ дает 56-мерное микровесовое представление группы типа $\E_7$, и разнообразные формы на этом представлении изучались во множестве работ, см., например, [30, 4]. Видимо, наиболее близкое к текущей работе задание этих форм встречается в [8]. Сходство, разумеется, не является полным~--- $y^{-\d}$ не является формой четвертой степени, а четырехлинейная форма из [8] не является симметричной. Однако оно слишком большое, чтобы являться случайным.

\begin{center}
{\bf Литература}
\end{center}

\begin{enumerate}

\item Борель А., {\it Свойства и линейные представления групп
Шевалле}, Семинар по алгебраическим группам, Мир, М., 1973, с.
9--59.

\item Бурбаки Н., {\it Группы и алгебры Ли}, Главы IV -- VI, Мир,
М., 1972.

\item Бурбаки Н., {\it Группы и алгебры Ли}, Главы VII -- VIII,
Мир, М., 1978.

\item Вавилов Н. А., Лузгарев А. Ю., {\it Нормализатор группы Шевалле типа $\E_7$}, Алгебра и анализ {\bf 27} (2015), №6, 57--88.

\item Вавилов Н. А., Лузгарев А. Ю., Певзнер И. М., {\it Группа
Шевалле типа $\E_6$ в $27$-мерном представлении}, Зап. научн.
семин. ПОМИ {\bf 338} (2006), 5--68.

\item Вавилов Н. А., Певзнер И. М., {\it Тройки длинных корневых
подгрупп}, Зап. научн. семин. ПОМИ {\bf 343} (2007), 54--83.

\item Вавилов Н. А., Семенов А. А., {\it Длинные корневые торы в
группах Шевалле}, Алгебра и анализ {\bf 24} (2012), №3, 22--83.

\item Лузгарев А. Ю., {\it Не зависящие от характеристики инварианты четвертой степени для $G(\E_7,R)$}, Вестник Санкт-Петербургского университета. Серия 1: математика, механика, астрономия (2013), №1, 43--50.

\item Лузгарев А. Ю., Певзнер И. М., {\it Некоторые факты из жизни
$\GL(5,\Z)$}, Зап. научн. семин. ПОМИ {\bf 305} (2003), 153--163.

\item О'Мира О., {\it Лекции о линейных группах}, Автоморфизмы
классических групп, Мир, М., 1976, с. 57--167.

\item О'Мира О., {\it Лекции о симплектических группах}, Мир, М.,
1979.

\item Певзнер И. М., {\it Геометрия корневых элементов в группах
типа $\E_6$}, Алгебра и анализ {\bf 23} (2011), №3, 261--309.

\item Певзнер И. М., {\it Ширина групп типа ${\E}_{6}$
относительно множества корневых элементов}, {\rm I}, Алгебра и
анализ {\bf 23} (2011), №5, 155--198.

\item Певзнер И. М., {\it Ширина групп типа ${\E}_{6}$
относительно множества корневых элементов}, {\rm II}, Зап. научн.
семин. ПОМИ {\bf 386} (2011), 242--264.

\item Певзнер И. М., {\it Ширина группы $\GL(6,K)$ относительно
множества квазикорневых элементов}, Зап. научн. семин. ПОМИ {\bf
423} (2014), 183--204.

\item Певзнер И. М., {\it Ширина экстраспециального унипотентного
радикала относительно множества корневых элементов}, Зап. научн.
семин. ПОМИ {\bf 435} (2015), 168--177.

\item Певзнер И. М., {\it Существование корневой подгруппы, которую данный элемент переводит в противоположную}, Зап. научн. семин. ПОМИ {\bf 460} (2017), 190--202.

\item Певзнер И. М., {\it Орбиты векторов некоторых представлений. {\rm I}}, Зап. научн. семин. ПОМИ {\bf 484} (2019), 149--164.

\item Певзнер И. М., {\it Орбиты векторов некоторых представлений. {\rm II}}, в печати.

\item Спрингер Т. А., {\it Линейные алгебраические группы},
Алгебраическая геометрия -- 4, Итоги науки и техн. Сер. Соврем.
проблемы мат. Фундам. направления {\bf 55}, ВИНИТИ, М., 1989, с.
5--136.

\item Стейнберг Р., {\it Лекции о группах Шевалле}, Мир, М., 1975.

\item Хамфри Дж., {\it Линейные алгебраические группы}, Наука, М.,
1980.

\item Хамфри Дж., {\it Введение в теорию алгебр Ли и их
представлений}, МЦНМО, М., 2003.

\item Aschbacher M., {\it The 27-dimensional module for $E_6$. {\rm I}}, Invent. Math {\bf 89} (1987), no. 1, 159--195.

\item Aschbacher M., {\it The 27-dimensional module for $E_6$. {\rm II}}, J. London Math. Soc {\bf 37} (1988), 275--293.

\item Aschbacher M., {\it The 27-dimensional module for $E_6$. {\rm III}}, Trans. Amer. Math. Soc. {\bf 321} (1990), 45--84.

\item Aschbacher M., {\it The 27-dimensional module for $E_6$. {\rm IV}}, J. Algebra {\bf 131} (1990), 23--39.

\item Aschbacher M., {\it Some multilinear forms with large isometry groups}, Geom.\ Dedicata {\bf 25} (1988), no. 1--3, 417--465.

\item Aschbacher M., {\it The geometry of trilinear forms}, Finite Geometries, Buildings and Related
topics, Oxford: Oxford Univ. Press (1990), 75--84.

\item Cooperstein B. N., {\it The fifty-six-dimensional module for $E_7$. I. A four form for $E_7$}, J. Algebra {\bf 173} (1995), no. 2, 361--389.

\item Krutelevich S., {\it Jordan algebras, exceptional groups, and Bhargava composition}, J. Algebra {\bf 314} (2007), no. 2, 924--977.

\item Springer T. A., {\it Linear algebraic groups}, Progress in
Mathematics {\bf 9}, Birkh\"auser Boston Inc., Boston, 1998.

\item Tits J., {\it Sur les constantes de structure et le th{\' e}or{\` e}me d'existence des alg{\` e}bres de Lie semi-simples}, Inst. Hautes {\' E}tudes Sci. Publ. Math. No. 31 (1966), 21--58.

\item Vavilov N. A., {\it A third look at weight diagrams}, Rend. Sem. Mat. Univ. Padova {\bf 104} (2000), 201--250.

\item Vavilov N. A., Plotkin E. B., {\it Chevalley groups over commutative rings. I. Elementary calculations}, Acta Applicandae Math. {\bf 45} (1996), 73--115.

\end{enumerate}

\end{document}